\documentclass[12pt]{amsart}

\usepackage{color}

\usepackage{graphicx,graphics}
\usepackage{fullpage,amssymb,amsfonts,amsmath,amstext,amsthm,amscd,enumitem,verbatim,tikz}
\usepackage[T1]{fontenc}
\setlist{itemsep=1mm}

\begin{document}
\newcommand{\note}[1]{\marginpar{\tiny #1}}

\newtheorem{theorem}{Theorem}
\theoremstyle{definition}
\newtheorem*{remark}{Remark}
\newtheorem{definition}{Definition}

\newcommand{\Z}{\mathbb{Z}}
\newcommand{\D}{\mathbb{D}}
\newcommand{\R}{\mathbb{R}}
\newcommand{\N}{\mathbb{N}}
\newcommand{\C}{\mathbb{C}}
\renewcommand{\H}{\mathbb{H}}
\renewcommand\Re{\operatorname{Re}}
\renewcommand\Im{\operatorname{Im}}
\newcommand{\diam}{\operatorname{diam}}
\newcommand{\dist}{\operatorname{dist}}
\newcommand{\ds}{\displaystyle}
%\numberwithin{equation}{section}
%\renewcommand{\theenumi}{(\roman{enumi})}
%\renewcommand{\labelenumi}{\theenumi}

\date{\today}
\title{Uniformly quasiregular maps of $\R^n$ \\ with Cantor Julia sets}
\author{Daniel A. Nicks}
\address{School of Mathematical Sciences, The University of Nottingham, Nottingham, NG7 2RD, United Kingdom}
%\textsc{\newline \indent {\includegraphics[width=1em,height=1em]{orcid2.png} {\normalfont https://orcid.org/0000-0002-9493-2970}}}}
\email{dan.nicks@nottingham.ac.uk}

\begin{abstract}
For any $n\ge 3$, we construct uniformly quasiregular self-maps of $\R^n$ of polynomial type for which the Julia set is a Cantor set.
\end{abstract}

\maketitle

\section{Introduction}

The iteration of uniformly quasiregular self-maps of $\overline{\R^n} = \R^n \cup \{\infty\}$ provides a natural generalisation of classical complex dynamics of rational maps on the Riemann sphere. Here a function $f\colon \overline{\R^n}\to \overline{\R^n}$ is called \emph{uniformly quasiregular} if there is a uniform bound on the distortion for the family of iterates of $f$; that is, there exists $K>1$ such that every iterate $f^k$ is $K$-quasiregular. The \emph{Fatou set} of $f$ is then, as usual, the set of points at which the iterates form a normal family, and the \emph{Julia set} $J(f)$ is its complement. Some further definitions are given in Section~\ref{sect:prelims}, while the reader is referred to \cite{Rickman} for additional  background regarding quasiregular mappings, and to \cite{HMM04} and the surveys \cite{Ber10}, \cite[Ch.~21]{IM02} and \cite[Ch.~4]{Sie04} for a fuller introduction to the dynamics of uniformly quasiregular maps.

The first examples of non-injective uniformly quasiregular self-maps of~$\overline{\R^n}$, with ${n\ge 3}$, were given by Iwaniec and Martin \cite{IwaMar96} via a conformal trap method. The breadth and complexity of the class of uniformly quasiregular maps is illustrated by Martin's short and elegant conformal trap construction \cite{Mar97} of a uniformly quasiregular map with an essentially arbitrary branch set $B$ --- here the \emph{branch set} is the set of points at which the map is not locally injective, and the only condition on $B$ is that it is the branch set of \emph{some} quasiregular map (not \emph{a priori} a uniformly quasiregular one). All the examples in \cite{IwaMar96} and \cite{Mar97} have Julia sets that are Cantor sets.
We also note that Martin and Peltonen \cite{MarPel10} proved that, on $\overline{\R^n}$, any quasiregular map can be written as a composition of a quasiconformal mapping and a uniformly quasiregular map with a Cantor Julia set.

A non-injective quasiregular map $f\colon \R^n\to \R^n$ is called \emph{polynomial type} if $|f(x)|\to \infty$ as $|x|\to\infty$, in which case we can extend it to a quasiregular self-map of $\overline{\R^n}$ by setting $f(\infty)=\infty$. The point at infinity is then a superattracting fixed point, and the Julia set equals the boundary of the basin of attraction. Mayer~\cite{May97} gave examples of polynomial type uniformly quasiregular maps  analogous to complex power mappings $z\mapsto z^d$. In particular, for each integer $d\ge 2$, there is a uniformly quasiregular `power map' $P_d\colon \R^n \to \R^n$ that satisfies
\begin{equation}\label{eqn:|P|}
|P_d(x)| = |x|^d.
\end{equation}
Note that \cite{May97} and \cite{May98} take $n=3$, but the same ideas works in all dimensions.  It follows immediately from \eqref{eqn:|P|} that $0$ and $\infty$ are superattracting fixed points of $P_d$ and that the Julia set is the unit sphere. The construction of $P_d$ is briefly outlined in Section~\ref{sect:proof}.  Mayer's method also gives uniformly quasiregular analogues of Chebyshev polynomials \cite{May98} and chaotic Latt\`es maps \cite{May97}. The former are polynomial type maps on $\R^n$ for which the Julia set is an $(n-1)$-dimensional ball, while the latter are self-maps of $\overline{\R^n}$ for which the Julia set is the whole space.

At this point observe that, for every $n\ge 2$, there are many examples of uniformly quasiregular self-maps of $\overline{\R^n}$ with Cantor Julia sets, while the polynomial type examples mentioned above have connected Julia sets. We ask the natural question: Are there uniformly quasiregular mappings $f\colon \R^n \to \R^n$ such that $J(f)$ is a Cantor set? Examples are known in low dimensions.  For $n=2$, there are of course many complex polynomials with Cantor Julia sets, including, for example, $z\mapsto z^d+c$ whenever $|c|$ is sufficiently large. We note also that the standard middle-thirds Cantor set can be realised as the Julia set of a uniformly quasiregular map of the plane, see for example~\cite[Theorem~1.3]{Fle19}. For $n=3$, Fletcher and Wu~\cite{FleWu15} exhibited uniformly quasiregular self-maps of $\R^3$ for which the Julia set is a wild Cantor set; see also~\cite{FleSto22, FSV25}. These constructions rely on an extension result due to Berstein and Edmonds that is only available in three dimensions. Pankka and Wu~\cite[Theorem~1.1.4]{PanWu19} establish the existence of similar wild Cantor sets in $n=4$ dimensions; their example is stated on~$\overline{\R^4}$, but it restricts to a self-map of~$\R^4$. 

The purpose of this paper is to give, for $n\ge 5$, the first examples of uniformly quasiregular mappings of~$\R^n$ with Cantor Julia sets. Moreover, we give a unified elementary construction for every dimension~$n\ge 3$.

\begin{theorem}\label{thm:main}
For each $n\ge 3$, there exists a hyperbolic uniformly quasiregular map ${f\colon \R^n\to \R^n}$ of polynomial type whose Julia set is a Cantor set.
\end{theorem}

Forthcoming work with Rose will produce examples with Julia sets of specified Hausdorff and packing dimensions.
 
 We remark that Fletcher and Stoertz~\cite{FleSto21} showed that if $f\colon \R^n\to\R^n$ is a uniformly quasiregular map of polynomial type with $J(f)$ a Cantor set, then every periodic point of $f$ must be repelling. The authors point out that this does not hold for Julia sets of rational maps on the Riemann sphere, which may contain parabolic fixed points.
 
An immediate application of Theorem~\ref{thm:main} is a generalisation of Fletcher and Vellis's result \cite[Theorem~1.2]{FleVel21} which states that if $E\subseteq\overline{\R^n}$ is a compact, uniformly perfect and uniformly disconnected set, then it is the Julia set of a hyperbolic uniformly quasiregular map ${f\colon\overline{\R^N}\to\overline{\R^N}}$, where $N = 2$ if $n = 2$ and $N = n + 1$ if $n \ge 3$. To see that this result holds with $\overline{\R^n}$  and $\overline{\R^N}$ replaced by $\R^n$ and $\R^N$, one uses in the proof the polynomial type maps from Theorem~\ref{thm:main} in place of the self-maps $G$ of $\overline{\R^N}$ given by Lemma~4.2 of \cite{FleVel21}. (Note that these maps do not have the property that $G^{-1}(\infty)=\infty$.)

Moreover, the proof given for the final assertion of Theorem~1.3 of \cite{FleVel21} appears to yield a hyperbolic uniformly quasiregular map of $\overline{\R^3}$, rather than of $\R^3$ as is claimed in the statement. To resolve this and obtain the claimed map of $\R^3$, one can similarly replace in the proof the (implicit) conjugation of a map $G$ of $\overline{\R^3}$ as above, with a map of $\R^3$ given by Theorem~\ref{thm:main}, noting that the Julia set of this map lies in a two-dimensional plane by construction.

\section{Preliminaries}\label{sect:prelims}

We will denote the open ball of radius $r$ centred at $a\in\R^n$ by $B(a,r)$ and, for $0<R<S$, we will denote annuli centred at $0$ by ${A(R,S) = \{x\in \R^n: R<|x|<S \}}$. The closure of a set $E\subseteq \R^n$ will be denoted by $\overline{E}$. A point $y\in\R^n$ will be taken to have co-ordinates $y=(y_1,\ldots, y_n)$.

Let $U \subset \R^n$ be a domain. A continuous map $f\colon U \rightarrow \R^n$ is called \emph{quasiregular} if it lies in the Sobolev space $W^{1}_{n, \textrm{loc}}(U)$ and if there exists $K_O \ge1$ such that
\[
|Df(x)|^n \le K_O J_{f}(x)
\]
for almost every $x\in U$. Here $|Df(x)|$ is the operator norm of the derivative $Df(x)$, and $J_{f}(x)$ denotes the Jacobian determinant.  If $f$ is quasiregular, then there also exists $K_I \ge 1$ such that
\[
J_{f}(x) \le K_I \inf_{|h|=1} |Df(x) (h)|^n
\]
for almost every $x\in U$. The smallest constants $K_O$ and $K_I$ for which the above hold are called the \emph{outer and inner dilatations} of $f$, denoted by $K_{O}(f)$ and $K_{I}(f)$. The \emph{dilatation} $K(f)$ of $f$ is the larger of $K_{O}(f)$ and $K_{I}(f)$, and we say that $f$ is $K$-quasiregular if $K(f)\le K$. An injective quasiregular map is called \emph{quasiconformal}.

In general, if $f$ is $K$-quasiregular on $\R^n$, then the iterates $f^k$ are also quasiregular, but the dilatation of the iterate $f^k$ may be as large as $K^k$. A map is said to be \emph{uniformly quasiregular} if there exists $K\ge 1$ such that every iterate is $K$-quasiregular with the same $K$.

\begin{definition}\label{defn:hyp} 
Let $f\colon \R^n\to\R^n$ be a non-injective uniformly quasiregular map. The \emph{branch set} $\mathcal{B}(f)$ is the set of points at which $f$ is not locally injective, and the \emph{post-branch set} is 
\[ \mathcal{P}(f) = \overline{\bigcup_{k \ge 0}  f^k(\mathcal{B}(f))  }. \]
The map $f$ is called \emph{hyperbolic} if the Julia set $J(f)$ is disjoint from $\mathcal{P}(f)$.
\end{definition}

The above definition is from \cite{FleVel21}, where it was shown that if the Julia set of a hyperbolic uniformly quasiregular map is totally disconnected, then it is uniformly disconnected. %\cite[Theorem~1.1]{FleVel21}.

We shall use the quasiconformal annulus theorem to interpolate between two given quasiconformal mappings. This result is due to Sullivan \cite{Sul79}, see also Tukia and V\"ais\"al\"a \cite{TukVai81}. We state a version where the domain is the difference of two quasiballs, which follows as in \cite{Fle19}. Here a \emph{quasiball} is the image of $B(0,1)$ under an ambient quasiconformal map of $\R^n$.

\begin{theorem}
\label{thm:ann}
Let $n\ge2$ and let $U$ and $V$ be %bounded 
quasiballs with $\overline{V} \subset U$, and let $f\colon W\to \R^n$ be a quasiconformal map on a neighbourhood $W$ of $\partial U \cup \partial V$ such that $f(W\cap U)$ lies in the bounded component of $\R^n\setminus f(\partial U)$ and also $f(W\setminus\overline{V})$ lies in the unbounded component of $\R^n\setminus f(\partial V)$. Then there exists a quasiconformal map $\widetilde{f}$ on a neighbourhood of $\overline U \setminus V$ that agrees with $f$ on a neighbourhood of $\partial U \cup \partial V$.
\end{theorem}

Informally, the mapping conditions above ensure that $f$ maps the inside of~$U$ to inside $f(\partial U)$ and maps the outside of~$\overline V$ to the outside of~$f(\partial V)$, so that the interpolation is possible topologically.

\section{Proof of Theorem~\ref{thm:main}}\label{sect:proof}

The construction begins by choosing a suitable Zorich map, so that the associated power map $P_2$ has certain convenient mapping properties. We remark that such choices make little difference to the function theory of the mappings, but can become important in a dynamical context as they determine the location and interaction of the branch set, its image, and any invariant hyperplanes. The choices made here may differ from those in other papers such as \cite{Ber10K, BFN25, May97, NicSix18, Tsa23, Tsa23a}. The function~$f$ that we seek will be obtained by modifying the power map~$P_2$ in a bounded region of $\R^n$ in order to create $2^{n-1}$ balls on each of which $f$ is an expanding similarity, while for any point outside these balls the orbit~$f^k(x)$ escapes to infinity.

We start with a bi-Lipschitz map $h$ from the $(n-1)$-cube $[0,1]^{n-1}$ onto the upper hemi-sphere $S_+^{n-1} = \{y\in\R^n : |y|=1, \, y_n \ge 0\}$. Note that $h$ must map the boundary of $[0,1]^{n-1}$ into $\{y\in S_+^{n-1} : y_n=0\}$ by the invariance of domain theorem. A Zorich map $\mathcal{Z}$ is then constructed as follows; see \cite{Zor67}, \cite[\S I.3.3]{Rickman} and \cite[\S 6.5.4]{IM02}. First, for $x=(x_1, \ldots, x_n) \in [0,1]^{n-1} \times \R$, we define
\begin{equation}\label{eqn:defnZ}
\mathcal{Z}(x) = e^{x_n}h(x_1, \ldots, x_{n-1}).
\end{equation}
This gives a homeomorphism of $[0,1]^{n-1} \times \R$ onto $\overline{\H^n_+}\setminus \{0\}$, where $\H^n_\pm = \{y\in \R^n : \pm y_n >0\}$ denotes an open half-space.  This mapping is then extended by repeatedly reflecting in the sides of the beam $[0,1]^{n-1} \times \R$ in the domain, and reflecting in the $y_n=0$ hyperplane in the codomain to obtain a quasiregular Zorich map $\mathcal{Z}\colon \R^n \to \R^n\setminus\{0\}$ (see the references cited earlier for further details). We note that if any of the first $n-1$ co-ordinates of $x$ is an integer, then $\mathcal{Z}(x)$ lies in the $y_n=0$ hyperplane.

For each integer $d\ge 2$, the uniformly quasiregular power map $P_d$  on $\R^n$ introduced by Mayer \cite{May97} is the solution to the Schr\"oder equation 
\begin{equation}\label{eqn:schroder}
\mathcal{Z}(dx) = P_d(\mathcal{Z}(x)),
\end{equation}
see also \cite[\S 21.4]{IM02}. 
Note that \eqref{eqn:|P|} follows from \eqref{eqn:defnZ} and \eqref{eqn:schroder}. Also, the degree of $P_d$ is $d^{n-1}$. Henceforth, we shall always take $d=2$.

Throughout the rest of this paper we impose the following additional condition on the map~$h$: for  $x\in [0,1]^{n-1}$ and  $j\in\{1,\ldots, n-1\}$, 
\begin{equation}\label{eqn:x_j=1/2}
x_j = \tfrac12 \ \ \implies \ \ h(x) \in \{y\in \R^n : y_j=0\}.
\end{equation} 
See the Remark at the end for an explicit example of such a map. Thus a hyperplane in $\R^n$ where one of the first $n-1$ co-ordinates is a half-integer is mapped by the Zorich map $\mathcal{Z}$ into the hyperplane where the corresponding co-ordinate is zero. In particular, $\mathcal{Z}(\frac12, \ldots, \frac12, 0) = (0, \ldots, 0, 1)$.  It follows that the power map $P_2$ quasiconformally maps each component of 
\[ \{ x\in\R^n : x_j\ne 0 \mbox{ for all } j=1,\ldots, n\} \]
 onto one of the half-spaces $\H^n_+$ or $\H^n_-$. To see this from \eqref{eqn:schroder}, note that each such component is the quasiconformal image under $\mathcal{Z}$ of $C\times \R$ where $C$ is an $(n-1)$-cube of side length $\frac12$.
 
 We now fix $R>4$ and consider the image under $P_2$ of certain `double quadrants' of the ball $B(0,R)$. Namely, for each $s=(s_1, \ldots, s_{n-1})$ in the index set $I=\{-1,1\}^{n-1}$, we let
 \[ Q_s = \{x \in \R^n : |x|<R, \ s_jx_j >0 \mbox{ for all } j=1, \ldots, n-1 \}. \]
 Then each set $Q_s$ is mapped quasiconformally by $P_2$ onto $B(0,R^2)\setminus E$, where $E$ is a subset of $B(0, R^2)\cap \{y_n=0\}$ that does not depend on the choice of $s\in I$. (In fact, 
 \begin{align*}
 E &= P_2\left(\{x \in \R^n : |x|<R, \ x_j =0 \mbox{ for some } j=1, \ldots, n-1 \}\right) \\
 &= \{ tz : 0< t < R^2, \ z\in h(N)\},
 \end{align*}
 where $N$ is the half of the boundary of $[0,1]^{n-1}$ consisting of points with $x_j=1$ for some $j=1, \ldots, n-1$.)
 
Let $a=(0, \ldots, 0, \frac12)\in \R^n$ and define a quasiconformal $\tilde{\phi}\colon B(0,1) \to B(0,1)$ that maps $B(a,\frac13)$ onto $B(0,\frac13)$ by translation, that is the identity on $A(\frac{9}{10}, 1)$, and interpolates between these on $B(0,\frac{9}{10}) \setminus B(a, \frac13)$ by the annulus theorem. Let $\phi\colon B(0,R^2) \to B(0,R^2)$ be given by $\phi(x) = R^2\tilde{\phi}(x/R^2)$. Then
\begin{enumerate} 
\item[(i)] $\phi(x)=x$ on $A\left(\tfrac{9}{10}R^2, R^2\right)$; and %\vspace{2mm}
\item[(ii)] $\phi(E) \subseteq \phi(\{ y_n=0\}) \subseteq A\left(\tfrac{1}{3}R^2, R^2\right)$.
\end{enumerate}

Returning to the sets $Q_s$, we will call $\{x\in\partial Q_s : |x|=R\}$ the \emph{spherical part} of the boundary $\partial Q_s$, and $\{ x\in \partial Q_s : x_j=0 \mbox{ for some } j=1, \ldots, n-1 \}$ the \emph{flat part} of $\partial Q_s$. For $\varepsilon>0$, we denote an $\varepsilon$-neighbourhood of $\partial Q_s$ in $Q_s$ as
\[ Q_s^\varepsilon = \{ x\in Q_s : \dist (x, \partial Q_s) < \varepsilon \}. \]
We claim that by continuity there exists  $0< \varepsilon < R/6\sqrt{n}$ such that, for each $s\in I$, the ring domain $Q_s^\varepsilon$ is mapped by $\phi \circ P_2$ into $A\left(\frac14 R^2, R^2\right)$. This is because $P_2$ maps the spherical part of $\partial Q_s$ onto $\partial B(0, R^2)$, near which $\phi$ is the identity, while $P_2$ maps the flat part of $\partial Q_s$ onto~$E$, that is then mapped by $\phi$ into $A\left(\frac13 R^2, R^2\right)$.
The claim follows.
Since $R>4$, we have that $B(0,R)$ is contained in the bounded component of $\R^n \setminus (\phi \circ P_2)(Q_s^\varepsilon)$. For each $s \in I$, take
\[ a_s = \left(\frac{s_1 R}{2\sqrt{n}}, \ldots, \frac{s_{n-1} R}{2\sqrt{n}}, 0 \right)\in Q_s \ \mbox{ and } \ B_s = B\left(a_s, \frac{R}{3\sqrt{n}}\right), \]
so that the ball $B_s$ lies in the bounded component of $\R^n \setminus Q_s^\varepsilon$. Now by the annulus theorem there exists a quasiconformal $f$ on $Q_s$ such that
\[ f(x) = \begin{cases} (\phi\circ P_2)(x), & \mbox{for } x\in Q_s^\varepsilon, \\
3\sqrt{n}(x- a_s), & \mbox{for } x\in B_s. \end{cases} \]
Moreover, for $|x|>R$ we set $f(x) = P_2(x)$. Then $f$ extends continuously to give a quasiregular map on $\R^n$, since near the flat part of each $\partial Q_s$ it is $\phi \circ P_2$ (for every $s\in I$) and near the sphere $\partial B(0,R)$ it is simply $P_2$. These regions do overlap, but $\phi \circ P_2 = P_2$ on the intersection.

The map $f$ is clearly of polynomial type with $\deg f = \deg P_2 = 2^{n-1}$. In fact, $f$ is uniformly quasiregular because
\begin{itemize}
\item $f$ maps $B_s$ onto $B(0,R)$ as a similarity, for each $s\in I$;
\item if $x\in B(0,R)\setminus \bigcup_{s\in I} B_s$, then $|f(x)|\ge R$; and
\item if $|x|\ge R$, then $f(x)=P_2(x)$ and thus $|f(x)|=|x|^2 > R$.
\end{itemize}
It follows that if an orbit $(f^k(x))$ ever leaves the set $\bigcup_{s\in I} B_s$ on which $f$ is conformal, then it spends at most one iteration in the intermediate zone $B(0,R)\setminus \bigcup_{s\in I} B_s$ before becoming trapped in $\R^n\setminus B(0,R)$, on which $f$ equals the uniformly quasiregular power map $P_2$.

Now let $X$ be the set of points $x$ for which the forward orbit $(f^k(x))$ is contained in $\bigcup_{s\in I} B_s$. Then
\[ X = \bigcap_{k\in\N} f^{-k}(\overline{B(0,R)}), \]
and thus $X$ is the Cantor set that is the attractor of the iterated function system of the $2^{n-1}$ contracting similarities which are the inverse branches of $f$ on $\overline{B(0,R)}$, given by
\[ y \mapsto \frac{y}{3\sqrt{n}} + a_s, \quad \mbox{ for } s\in I. \]
Note that for $x\notin X$ we have $f^k(x)\to\infty$ as $k\to\infty$. Therefore, $X$ is the Julia set of $f$.

It only remains to show that $f$ is hyperbolic as in Definition~\ref{defn:hyp}. The branch set is disjoint from $\bigcup_{s\in I} B_s$ as $f$ is injective on each ball $B_s$, and it follows from the bullet points above that the post-branch set $\mathcal{P}(f)$ is also disjoint from $\bigcup_{s\in I} B_s$. Therefore, since the Julia set $J(f)= X \subseteq \bigcup_{s\in I} B_s$, we conclude that $f$ is indeed hyperbolic.

\begin{remark}
The following standard example of a bi-Lipschitz map $h\colon [0,1]^{n-1} \to S_+^{n-1}$ satisfies the additional condition \eqref{eqn:x_j=1/2}. For $x\in [0,1]^{n-1}$, let $x'=x-(\frac12, \ldots, \frac12)$, recall that $|x'|=\sqrt{x'_1{}^2+\ldots+x'_{n-1}{}^2}$
 and let now $\| x' \|_\infty = \max\{|x'_1|, \ldots, |x'_{n-1}|\}$. Then we can take
\[ h(x) = \left(\frac{x'}{|x'|}\sin\left(\pi \|x' \|_\infty\right), \cos\left(\pi \|x' \|_\infty\right) \right) \in \R^{n-1}\times \R. \]
\end{remark}

\subsection*{Acknowledgments}
The author thanks Alastair Fletcher and Peter Rose for interesting and useful discussions that prompted this work.

\bibliography{book}
\bibliographystyle{plain}
\end{document}